\documentclass[12TP,draft]{article}

\begin{document}

\begin{center}
\LARGE\noindent\textbf{ On Cycles through Vertices of Large Semidegree in Digraphs }\\

\end{center}
\begin{center}
\noindent\textbf{S.Kh. Darbinyan and I.A. Karapetyan}\\

Institute for Informatics and Automation Problems, Armenian National Academy of Sciences

E-mails: samdarbin@ipia.sci.am, isko@ipia.sci.am\\
\end{center}

\textbf{Abstract}\\

 Let $D$ be a strong digraph on $n=2m+1\geq 5$ vertices. In this paper we show that if $D$ contains a cycle of length $n-1$, then $D$  has  also a cycle which contains all vertices with in-degree and out-degree at least $m$ (unless some extremal cases).


\noindent\textbf{1. Introduction }\\

The digraph $D$ is hamiltonian  if it contains a hamiltonian cycle, i.e. a cycle of length $|V(D)|$. A set $S$ of vertices in a digraph $D$ (an undirected graph $G$) is said to be cyclable in $D$ (in $G$) if $D$ ($G$) contains a cycle through all vertices of $S$.

There are many well-known conditions which guarantee the cyclability of a set of vertices in undirected graph. Most of them can be seen as restrictions of hamiltonian conditions to the considered set of vertices (See \cite{[4], [5], [15], [16], [18]}). However, for general digraphs, relatively few degree conditions are known to guarantee hamiltonisity in digraphs (See \cite {[2], [3], [7], [9], [13], [14], [17], [19]}). The more general and classical ones is the following theorem of M. Meyniel:  

\noindent\textbf{Theorem A} \cite{[13]}. If $D$ is a strong digraph of order $n\geq 2$ and $d(x)+d(y)\geq 2n-1$ for all pairs of nonadjacent vertices in $D$, then $D$ is hamiltonian .\\

In \cite{[8]} the first author proved the following:

\noindent\textbf{Theorem B} \cite{[8]}. Let $D$ be a strong digraph of order $n\geq 3$. If $d(x)+d(y)\geq 2n-1$ for any two non-adjacent vertices $x,y\in V(D)-\{z_0\}$, where $z_0$ is  some vertex of $D$, then $D$ is hamiltonian or contains a cycle of length $n-1$.\\
The following result is immediately corollary of Theorem B.

\noindent\textbf{Corollary} \cite{[8]}. Let $D$ be a strong digraph of order $n\geq 3$. If $D$ has $n-1$ vertices of degree at least $n$, then $D$ is a hamiltonian or contains a cycle of length $n-1$.\\

A Meyniel set $M$ is a subset of $V(D)$ such that  $d(x)+d(y)\geq 2n-1$ for every pair of vertices  $x$, $y$ in $M$ which are nonadjacent in $D$. In \cite{[4]}, K. A. Berman and X. Liu improved Theorem B proving the following generalization of well- known Meyniel's theorem.\\

  \noindent\textbf{Theorem C} \cite{[4]}. Let $D$ be a digraph of order $n$. If $D$ is strongly connected, then every Meyniel set $M$ lies in a cycle.\\
Theorem C also generalizes the classical theorems A. Ghouila-Houri \cite{[11]} and D.R. Woodall \cite{[19]}.\\

The digraph $D$ is $S$-strongly connected if for any pair $x,y$ of distinct vertices of $S$ there exists a path from $x$ to $y$ and a path from $y$ to $x$ in $D$ (See \cite{[12]}). H. Li, E. Flandrin and J. Shu \cite{[12]} proved the following generalization of Theorem C.

\noindent\textbf{Theorem D} \cite{[12]}. Let $D$ be a digraph of order $n$ and $M$ be a Meyniel set in $D$. If $D$ is $M$-strongly connected, then $D$ contains a cycle through all vertices of  $M$. \\

C. Thomassen \cite{[17]} (for $n=2k+1$) and first author \cite{[7]} (for $n=2k$) proved the following: 

\noindent\textbf{Theorem E} \cite{[17], [7]}. If $D$ is a digraph of order $n\geq 5$ with minimum degree at least $n-1$ and with minimum semi-degree at least $n/2-1$, then $D$ is hamiltonian (unless some extremal cases which are characterized). \\

 We put as a question to known if this result of C. Thomassen and first author has a cyclable version.

Let $D$ be a digraph of order  $n=2m+1$. A  Thomassen set $T$ is a subset of $V(D)$ such that $d^+(x)\geq m$ and $d^-(x)\geq m$ for every $x\in T$, we denote the vertices of $T$ by $T$-vertices. The cycle containing all vertices of $T$ is called an $T$-cycle.

In this paper we prove the following two theorems which provide some support for the above question. 

\noindent\textbf{Theorem 1}. Let $D$ be a 2-strong digraph of order $n=2m+1\geq 3$. Then any two $T$-vertices $x$ and $y$ are on a common cycle in $D$.\\

\noindent\textbf{Theorem 2}. Let $D$ be a strong digraph of order $n=2m+1\geq 3$. If $D$ contains a cycle of length $n-1$,  then  $D$ also contains a cycle containing all vertices with in-degree and out-degree at least $m$ unless some extremal cases.

 Our proofs are based on the arguments of \cite{[17], [7]}.\\

\noindent\textbf{2.Terminology and notations}\\

We shall assume that the reader is familiar with the standard terminology on directed graphs (digraphs) and refer the reader to monograph of Bang-Jensen and Gutin \cite{[1]} for terminology not discussed here. In this paper we consider finite digraphs without loops and multiple arcs. For a digraph $D$, we denote by $V(D)$ the vertex set of $D$ and by  $A(D)$ the set of arcs in $D$. The order $|D|$ of $D$ is the number of its vertices. Often we will write $D$ instead of $A(D)$ and $V(D)$. The arc of a digraph $D$ directed from $x$ to $y$ is denoted by $xy$. For disjoint subsets $A$ and  $B$ of $V(D)$  we define $A(A\rightarrow B)$ \, as the set $\{xy\in A(D) / x\in A, y\in B\}$ and $A(A,B)=A(A\rightarrow B)\cup A(B\rightarrow A)$. If $x\in V(D)$ and $A=\{x\}$ we write $x$ instead of $\{x\}$. If $A$ and $B$ are two disjoint subsets of $V(D)$ such that every vertex of $A$ dominates every vertex of $B$, then we say that $A$ dominates $B$, denoted by $A\rightarrow B$. The out-neighborhood of a vertex $x$ is the set $N^+(x)=\{y\in V(D) / xy\in A(D)\}$ and $N^-(x)=\{y\in V(D) / yx\in A(D)\}$ is the in-neighborhood of $x$. Similarly, if $A\subseteq V(D)$ then $N^+(x,A)=\{y\in A / xy\in A(D)\}$ and $N^-(x,A)=\{y\in A / yx\in A(D)\}$. We call the vertices in $N^+(x)$, $N^-(x)$, the out-neighbors and in-neighbors of $x$. The out-degree of $x$ is $d^+(x)=|N^+(x)|$ and $d^-(x)=|N^-(x)|$ is the in-degree of $x$. The out-degree and in-degree of $x$ we call its semi-degrees. Similarly, $d^+(x,A)=|N^+(x,A)|$ and $d^-(x,A)=|N^-(x,A)|$. The degree of the vertex $x$ in $D$ defined as $d(x)=d^+(x)+d^-(x)$ (similarly, $d(x,A)=d^+(x,A)+d^-(x,A)$). The subdigraph of $D$ induced by a subset $A$ of $V(D)$ is denoted by $\langle A\rangle$. The path (respectively, the cycle) consisting of the distinct vertices $x_1,x_2,\ldots ,x_m$ ( $m\geq 2 $) and the arcs $x_ix_{i+1}$, $i\in [1,m-1]$  (respectively, $x_ix_{i+1}$, $i\in [1,m-1]$, and $x_mx_1$), is denoted  $x_1x_2\cdots x_m$ (respectively, $x_1x_2\cdots x_mx_1$). For a cycle  $C_k=x_1x_2\cdots x_kx_1$, the subscripts considered modulo $k$, i.e. $x_i=x_s$ for every $s$ and $i$ such that  $i\equiv s\, (\hbox {mod} \,k)$. If $P$ is a path containing a subpath from $x$ to $y$ we let $P[x,y]$ denote that subpath. Similarly, if $C$ is a cycle containing vertices $x$ and $y$, $C[x,y]$ denotes the subpath of $C$ from $x$ to $y$.
 A digraph $D$ is strongly connected (or just strong) if there exists a path from $x$ to $y$ and a path from $y$ to $x$ in $D$ for every choice of distinct vertices $x$,\,$y$ of $D$. A digraph $D$ is $k$-connected, $k>0$ (or $k$-strong) if $|V(D)|\geq k+1$ and deletion of fewer than $k$ vertices always results in a strong digraph. For an undirected graph $G$, we denote by $G^*$ symmetric digraph obtained from $G$ by replacing every edge $xy$ with the pair $xy$, $yx$ of arcs. $K_n$ (respectively, $K_{p,q}$)   denotes the complete graph of order $n$ (respectively, complete bipartite graph  with partite sets of cardinalities $p$ and $q$), and   $\overline K_{n}$ denotes the complement of complete undirected graph of order $n$. Two distinct vertices $x$ and $y$ are adjacent if $xy\in A(D)$ or $yx\in A(D) $ (or both). We denote by $a(x,y)$ the number of arcs between the vertices $x$ and $y$. In particular, $a(x,y)=0$ (respectively, $a(x,y)\not=0$) means that $x$ and $y$ are not adjacent (respectively, are adjacent). 

For integers $a$ and $b$, $a\leq b$, let $[a,b]$  denote the set of all integers which are not less than $a$ and are not greater than $b$.

\noindent\textbf{3. Preliminaries }\\

The following well-known simple lemmas is the basis of our results and other theorems on directed cycles and paths in digraphs. It we  will be used extensively in the proofs of our results.\\

\noindent\textbf{Lemma 1} \cite{[10]}. Let $D$ be a digraph on $n\geq 3$ vertices containing a cycle $C_m$, $m\in [2,n-1] $. Let $x$ be a vertex not contained in this cycle. If $d(x,C_m)\geq m+1$, then  $D$ contains a cycle $C_k$ for all  $k\in [2,m+1]$.   \\

\noindent\textbf{Lemma 2} \cite{[6]}. Let $D$ be a digraph on $n\geq 3$ vertices containing a path $P:=x_1x_2\ldots x_m$, $m\in [2,n-1]$ and let $x$ be a vertex not contained in this path. If one of the following conditions holds:

 (i) $d(x,P)\geq m+2$; 

 (ii) $d(x,P)\geq m+1$ and $xx_1\notin D$ or $x_mx_1\notin D$; 

 (iii) $d(x,P)\geq m$, $xx_1\notin D$ and $x_mx\notin D$;

\noindent\textbf{}then there is an  $i\in [1,m-1]$ such that $x_ix,xx_{i+1}\in D$, i.e., $D$ contains a path $x_1x_2\ldots x_ixx_{i+1}\ldots x_m$ of length $m$  (we say that  $x$ can be inserted into $P$ or the path $x_1x_2\ldots x_ixx_{i+1}\ldots x_m$ is extended from $P$ with $x$ ). \\

\noindent\textbf{4. Main results}\\

\noindent\textbf{Theorem 1}. Let $D$ be a 2-strong digraph of order $n=2m+1\geq 3$. Then any two $T$-vertices $x$ and $y$ are on a common cycle in $D$.

\noindent\textbf{Proof}. Suppose, on the contrary, that there are two $T$-vertices  $x$ and $y$ which are not on common cycle. The vertices $x$ and $y$ are not adjacent, otherwise, if for example there is the arc $xy$, then using a path from $y$ to $x$ that necessarily exists from strong property of $D$, we get a contradiction. Denote $R:=N^+(x)\cap N^-(y)$ and $Q:=N^+(y)\cap N^-(x)$. The  assumption that $x$ and $y$ are $T$-vertices implies that $Q$ and $R$ (both) are nonempty. If $R\not= Q$ or $|R|\geq 2$, then the theorem is true. Assume that $R=Q=\{z\}$. Then  $V(D)=A\cup B\cup \{x,y,z\}$, where $A:=N^+(x)\setminus \{z\}$ and $B:=N^-(y)\setminus \{z\}$. Let the sets $A$ and $B$ (both) are not empty, i.e., $n\geq 5$. It is easy to see that $A(A\rightarrow B)=\emptyset$. In particular, $D$ is not 2-strong which is a contradiction. \fbox \\\\

For the next theorem we need the following definitions.

  \noindent\textbf{Definition 1}. $D_7$ is a digraph (see \cite{[1], [17]}) with vertex set $V(D_7)=\{x_1,x_2,x_3,x_4,x_5,x,y\}$ such that $N^+(x_1)=\{x_2,x_5,y\}$, $N^+(x_2)=\{x_3,x_4,y\}$, $N^+(x_3)=\{x_2,x_4,x\}$, $N^+(x_4)=\{x_3,x_5,x\}$, $N^+(x_5)=\{x_1,x,y\}$, $N^+(x)=\{x_1,x_2,x_3\}$ and $N^+(y)=\{x_1,x_4,x_5\}$.

\noindent\textbf{Definition 2}. $D_5$ is a digraph (see \cite{[1], [17]}) with vertex set $V(D_5)=\{x_1,x_2,x_3,x,y\}$ such that $N^+(x_1)=\{x_2,y\}$, $N^+(x_2)=\{x_3,x\}$, $N^+(x_3)=\{x,y\}$,  $N^+(x)=\{x_1,x_2\}$ and $N^+(y)=\{x_1,x_3\}$.

We denote by $L_1$ the set of three digraphs obtaining from $D_5$ by adding the arc $x_1x_3$ or $x_3x_1$ (or both).

\noindent\textbf{Definition 3}. By $L_2$ we denote the set of digraphs $D$ with vertex set $V(D)=\{x_1,x_2,\ldots ,x_{2m},x\}$ and with the following properties:

i. $D$ contains a cycle $x_1x_2\ldots x_{2m}x_1$ of length $2m$ and the vertices $x$ and $x_{2m}$ are not adjacent;

ii.  $N^+(x)= N^+(x_{2m}) =\{x_1,x_2, \ldots ,x_{m}\}$ and $N^-(x)= N^-(x_{2m}) =\{x_{m},x_{m+1}, \ldots ,x_{2m-1}\}$;

iii. $A(\{x_1,x_2, \ldots ,x_{m-1}\}\rightarrow \{x_{m+1},x_{m+2}, \ldots ,x_{2m-1}\})=\emptyset$, the induced subdigraphs $\langle \{x_1,x_2, \ldots ,$ $x_{m}\}\rangle$ and $\langle \{x_{m},x_{m+1}, \ldots ,x_{2m-1}\}\rangle$ are arbitrary and one may add any number of arcs that go from $\{x_{m+1},x_{m+2}, \ldots ,x_{2m-1}\}$ to $\{x_1,x_2, \ldots , x_{m}\}$. (Note that the digraphs from $L_2$ is not 2-strong and $x$, $x_{2m}$ are $T$-vertices which are not in common cycle.\\

In further, by $H$ we denote a hamiltonian cycle in $D$.\\

\noindent\textbf{Theorem 2}. Let $D$ be a strong digraph of order $n=2m+1\geq 3$ and $D$ contains a cycle of length $n-1$.  Then one of the following holds:

i. $D$ contains a cycle containing all vertices with in-degree and out-degree at least $m$;

ii. $D$ is isomorphic to digraphs $D_5$ or $D_7$ or belongs to the set $L_1\cup L_2$;

iii. $K^*_{m,m+1}\subseteq D \subseteq [K_m+\overline K_{m+1}]^*$;

iv. $D$ contains a cycle $C:=x_1x_2\ldots x_{2m}x_1$ of length $n-1$, and if $x\notin V(C)$ and $x$ is not adjacent with the vertices $x_{l_1}, x_{l_2},\ldots ,x_{l_j}$, $j\geq 3$, then $x_{l_i-1}x, xx_{l_i+1}\in D$ and $N^+(x)=N^+(x_{l_i})$ and $ N^-(x)=N^-(x_{l_i})$ for all $i\in [1,j]$. In particular, $\{x_{l_1}, x_{l_2},\ldots ,x_{l_j},x \}$ is an independent set of vertices.

 \noindent\textbf{Proof}. The proof is by contradiction. Suppose that Theorem 2 is false, in particular, $D$ is not hamiltonian. Let $C:=x_1x_2\ldots x_{n-1}x_1$  be an arbitrary cycle of length $n-1$ in $D$ and let the vertex $x$ is not containing this cycle $C$. Then $x$ is a $T$-vertex. Since $C$ is a longest cycle, using Lemmas 1 and 2, we obtain the following claim:\\

\noindent\textbf {Claim 1.} (i). $d(x)=n-1$ and there is a vertex $x_l$, $l\in [1,n-1]$ which is not adjacent with $x$.

 (ii). If $x_ix\notin D$, then $xx_{i+1}\in D$ and if $xx_i\notin D$, then $x_{i-1}x\in D$, where $i\in [1,n-1]$.

(iii). If the vertices $x$ and $x_i$ are not adjacent, then $x_{i-1}x,\,xx_{i+1}\in D$ and $d(x_i)=n-1$. \fbox \\\\

 By Claim 1(i), without loss of generality, we may assume that the vertices $x$ and $x_{n-1}$ are not adjacent. For convenience, let $p:=n-2$ and $y:=x_{n-1}$. We have $ yx_1,$ $x_py\in D$ and  $x_px$, $xx_1\in D$ by Claim 1(iii). Therefore $y$ is a $T$-vertex and $d(y)=n-1$.\\

\noindent\textbf{Claim 2}. At least two vertices of $C$ are not adjacent with $x$ unless $D$ is isomorphic to $D_5$ or $D_7$ or belongs to the set $L_1\cup L_2$.

\noindent\textbf{Proof}. We prove Claim 2 by contradiction. Let $C:=x_1x_2 \ldots x_{n-1}x_1$. Then, by Lemma 1, $d(x)=n-1$ and $d^+(x)=d^-(x)=m$ since $D$ is not hamiltonian. It is easy to see that some vertex $x_i$ (say, $y:=x_{n-1}$) is not adjacent with $x$. Then, by Claim 1(iii), $x_{p}x$, $xx_{1}\in D$. If $y$ is not a $T$-vertex, then the cycle $x_1x_2 \ldots x_{n-2}yx_1$ contains all $T$-vertices. So, we can assume that $y$ is a $T$-vertex. Then $d(y)=n-1$ (by Lemma 1) and $d^+(y)=d^-(y)=m$. From our assumption it follows that
$$
N^+(x)=\{x_1,x_2, \ldots , x_m\} \quad \hbox{and} \quad N^-(x)=\{x_m,x_{m+1}, \ldots , x_{p}\}. \eqno (1)
$$

We first prove that there is a vertex $x_k$, $k\in [2,p-1]$, which is not adjacent with $y$. Assume that it is not the case.  Then 
$$
N^+(y)=\{x_1,x_2, \ldots , x_m\} \quad \hbox{and} \quad N^-(y)=\{x_m,x_{m+1}, \ldots , x_{p}\}. \eqno (2)
$$
Since $D$ is not hamiltonian we have
$$
A(\{x_1,\ldots , x_{m-1}\}\rightarrow \{x_{m+1},\ldots , x_{p}\})=\emptyset, \eqno (3)
$$
for otherwise, if $x_ix_j\in D$, where $i\in [1,m-1]$ and $j\in [m+1,n-2]$, then by (1) and (2), $H=x_1\ldots x_ix_j \ldots x_{p}xx_{i+1}\ldots $ $ x_{j-1}yx_1$ is a hamiltonian cycle. Therefore 
$$
A(\{x_1,\ldots , x_{m-1},x,y\}\rightarrow \{x_{m+1},\ldots , x_{n-2}\})=\emptyset, 
$$
i.e., $D$ belongs to the set $L_2$ which is a contradiction. 

Thus  there is a vertex $x_k$ with $k\in [2,p-1]$ which is not adjacent with $y$. By Claim 1(iii), $x_{k-1}y, yx_{k+1}\in D$. Observe that $x_k$ also is a $T$-vertex. If $k\in [m+1,p-1]$, then $m\geq 3$ and from $d^-(x_k,\{x,y\})=0$ it follows that there is a vertex $x_i$, $i\in [1,m-1]$, such that $x_ix_k\in D$. Therefore $H=x_1\ldots x_ix_k\ldots x_{p}xx_{i+1}\ldots $ $ x_{k-1}yx_1$, a contradiction. So, we can assume that $k\leq m$. Similarly, we can assume that $k\geq m$. Therefore remains to consider the case when $m=k$ and the vertex $y$ is adjacent with all vertices  of $P\setminus \{x_m\}$. If $n=5$, i.e., $m=2$, then  $x_1y,yx_3\in D$ and $x_2x_1\notin D$, $x_3x_2\notin D$, i.e., $D$ isomorphic to well-known digraph $D_5$ or $D\in L$, since  if we add the arc $x_1x_3$ or $x_3x_1$ (or both) to $D_5$, then the resulting digraph also is not hamiltonian, i.e., $D\in L_1$. Assume that $m\geq 3$. It is not difficult to see that 
$$
d(x_m, \{x_1, x_{p}\})=0 \,\,\, \hbox{and}\,\,\,A(\{x_1,\ldots , x_{m-2}\}\rightarrow x_m)=A(x_m\rightarrow  \{x_{m+2},\ldots , x_{p}\})=\emptyset, \eqno (4)
$$
in particular, $x_m$ is not adjacent with $x_1$ and $x_{p}$. Therefore
$$
\{x_{m+1},\ldots , x_{p-1}\}\rightarrow x_m \rightarrow  \{x_{2},\ldots , x_{m-1}\}. \eqno (5)
$$
This implies that $x_{p}$ and $x_1$ are $T$-vertices since $x_1\ldots x_{m-1}yx_{m+1}\ldots x_{p-1}x_mxx_1$ (respectively,  $x_2\ldots $ $ x_{m-1}y$ $x_{m+1}\ldots x_{p}xx_mx_2$) is a cycle of length $n-1$ which does not contain $x_{p}$ (respectively, $x_1$).

Now we consider the vertex $y$. If $x_{p-1}y \in D$, then $xx_{p}\notin D$ and $yx_{p}\notin D$ imply that $x_ix_{p}\in D$ for some $i\in [1,m-1]$, and hence $H=x_1\ldots x_ix_{p}xx_{i+1}\ldots x_{p-1}yx_1$, a contradiction. So, we can assume that $x_{p-1}y \notin D$ and, similarly, $yx_2\notin D$, i.e., $yx_{p-1},x_2y\in D$. Using Lemma 2 we obtain that 
$$
\{x_{1},x_2, \ldots , x_{m-1}\}\rightarrow y \rightarrow  \{x_{m+1},x_{m+2},\ldots , x_{p}\}. \eqno (6)
$$
It is not difficult to see that $d^+(x_1, P[x_3,x_{m+1}])=0$, for otherwise, if $x_1x_i\in D$, $i\in [3,m]$, then by (1) and (6),  $H=x_1x_i\ldots x_{p}xx_2\ldots x_{i-1}yx_1$, and  if $x_1x_{m+1}\in D$, then by (1), (5) and (6),   $H=x_1x_{m+1}\ldots x_{p}xx_mx_2\ldots x_{m-1}yx_1$, which is a contradiction. Similarly, we can show that $d^-(x_{p},$ $P[x_{m-1},$ $x_{p-2}])=0$. Therefore 
$$
N^+(x_1)=\{x_2,y,x_{m+2},x_{m+3}, \ldots , x_{p}\} \quad \hbox{and} \quad N^-(x_{p})=\{x_{p-1},y,x_1,x_{2}, \ldots , x_{m-2}\}. \eqno (7)
$$
 By (7), (5) and (6) it is easy to see that  $x_1\ldots x_{m-2}x_{p}yx_{m+1}\ldots x_{p-1}x_mxx_1$ (respectively, $x_1 x_{m+2}$ $\ldots x_{p}xx_{m}$ $x_2 \ldots x_{m-1}yx_1$) is a cycle of length $n-1$, which does not contain $x_{m-1}$ (respectively, $x_{m+1}$). This means that $x_{m-1}$ and $x_{m+1}$ are $T$-vertices.

Now we will consider the vertex $x_{m-1}$. Then $x_{m-1}x_i\notin D$ for all $i\in [m+2,p]$ (for otherwise, by (5), $H=x_1\ldots x_{m-1}x_i\ldots x_{p}yx_{m+1}\ldots x_{i-1}x_mxx_1$) and $x_{m-1}x_1\notin D$ (for otherwise, $H=x_1\ldots x_{m-2}yx_{m+1}$ $\ldots x_{p}xx_mx_{m-1}x_1$ by (5) and (6)). Thus we have $d^+(x_{m-1}, \{x_1,x,x_{m+2},\ldots ,x_{p}\})$ $=0$. Therefore   
$$
x_{m-1}\rightarrow \{x_2,\ldots ,x_{m-2},y,x_{m},x_{m+1}\}. \eqno (8)
$$
Now, if $m\geq 4$, then by (7), (1), (8) and (5) we have $H=x_1x_{p}xx_{m-1}x_{m+1}\ldots x_{p-1}x_mx_2\ldots x_{m-2}yx_1$, which is a contradiction.

Therefore $m=3$, i.e., $n=7$. From (4), (5) and (7) we obtain that $x_4x_3, x_3x_2, x_1x_5\in D$, $x_1$ and $x_5$ are $T$-vertices and $d(x_3,\{x_1, x_5\})=0$. It is easy to see that $d^+(x_2,\{x_1, x_5\})=d^+(x_5,\{x_2, x_4\})=0$. From this we conclude that $x_5x_1\in D$. Now we see that $x_1x_5yx_4x_3xx_1$ is a cycle of length $n-1$ which does not contain $x_2$. This means that $x_2$ is a $T$-vertex and $d^+(x_2)=d^-(x_2)=3$. Since $d^+(x_2,\{x, x_1, x_5\})=0$, it follows that $x_2x_4\in D$. Therefore $D$ is isomorphic to digraph $D_7$. Claim 2 is proved. \fbox \\\\

\noindent\textbf{Claim 3}. Let $x_{p-1}x, yx_p\in D$ and for some $k\in [2,p-2]$ $x_k$ and $y$ are not adjacent. Then $x_k $ and $x_p$ also are not adjacent.

\noindent\textbf{Proof}. Since $x_k$ and $y$ are not adjacent it follows that $x_{k-1}y, yx_{k+1}\in D$ (by Claim 1(iii)). Now if $x_kx_p\in D$, then $H=x_1\ldots x_kx_pyx_{k+1}\ldots x_{p-1}xx_1$; and if $x_px_k\in D$, then $H=x_1\ldots x_{k-1}yx_px_k\ldots x_{p-1}xx_1$. In each case we have obtained a hamiltonian cycle, which is a contradiction.                               \fbox \\\\

 \noindent\textbf{Claim 4}. If $x_{p-1}x$ and $yx_p \in D$, then $d(x_i,\{x,y\})\geq 1$ for all $i\in [2,p-2]$.

\noindent\textbf{Proof}. Suppose, on the contrary, that $d(x_i,\{x,y\})=0$ for some $i\in [2,p-2]$. Then by Claim 1(iii), $x_{i-1}\rightarrow \{x,y\}\rightarrow x_{i+1}$, and by Claim 3 the vertices $x_i$ and $x_p$ are not adjacent. Now, since $x_i$ is a $T$-vertex and cannot be inserted into  $P[x_1,x_{i-1}]$ and into $P[x_{i+1},x_{p-1}]$, using Lemma 2 we obtain that
$$ 
p+1=d(x_i)=d(x_i,P[x_1,x_{i-1}])+d(x_i,P[x_{i+1},x_{p-1}])\leq i+p-i=p,
$$
a contradiction.  \fbox \\\\

 \noindent\textbf{Claim 5}. If $x_{p-1}x\in D$, then the vertices $y$ and $x_{p-1}$ are adjacent.

\noindent\textbf{Proof}. Suppose, on the contrary, that $y$ and $x_{p-1}$ are not adjacent. Then by Claim 1(iii), $x_{p-2}y, yx_p\in D$. If $x_px_{p-1}\in D$, then $H=x_1\ldots x_{p-2}yx_px_{p-1}xx_1$, a contradiction. So, we can assume that $x_px_{p-1}\notin D$. Moreover, if $xx_i\in D$ with $i\in [2,p-2]$, then $x_{i-1}x_{p-1}\notin D$ (for otherwise, we would be have a hamiltonian cycle $H=x_1\ldots $ $x_{i-1}x_{p-1}x_pxx_i\ldots x_{p-2}yx_1$). Recall (by Claim 2) that there is a vertex $x_l$ with $l\in [2,p-2]$ which is not adjacent with $x$. Note that $x_{l-1}x$ and $xx_{l+1}\in D$ by Claim 1(iii). Since $x$ is a $T$-vertex, it follows that $d^+(x,P[x_2,x_{p-2}])\geq m-2$. If we consider the vertex $x_{p-1}$, then from  $d^-(x_{p-1},\{y,x_{p}\})=0$ and the above observation it follows that 
$$xx_{p-1}\,\,\, \hbox {and}\,\,\, x_{l-1}x_{p-1}\in D.   \eqno (9) $$ 
Hence $x_px_l\notin D$ (for otherwise, if $x_px_l\in D$, then $H=x_1\ldots x_{l-1}xx_{p-1}x_px_l\ldots x_{p-2}yx_1$).  We consider the following two cases.

\noindent\textbf{Case 5.1}. $l\leq p-3$. Then it is not difficult to see that the vertices $x_l$ and $x_{p-1}$ are not adjacent. Indeed, if $x_{p-1}x_l\in D$, then $H=x_1\ldots x_{l-1}x_{p-1}x_l\ldots x_{p-2}yx_pxx_1$ by (9); and if $x_lx_{p-1}\in D$, then $H=x_1\ldots x_lx_{p-1}x_pxx_{l+1}\ldots x_{p-2}yx_1$, which is a contradiction. From this we have
$$
p+1=d(x_l)=d(x_l,P[x_1,x_{l-1}])+d(x_l,P[x_{l+1},x_{p-2}])+d(x_l,\{y,x_p\}).   \eqno (10)
$$

Now we show that
$$ x_lx_p \,\, \hbox{and} \,\,  x_{p-2}x_l\in D. \eqno (11)$$

Let first $yx_l\in D$. Then $x_lx_1\notin D$ (for otherwise, $H=x_1\ldots x_{l-1}xx_{l+1}\ldots x_pyx_lx_1$ is a hamiltonian cycle, a contradiction). Since the vertex $x_l$ cannot be inserted into $P[x_1,x_{l-1}]$ and $P[x_{l+1},x_{p-2}]$, from (10), $x_px_l\notin D$ and Lemma 2 it follows that $d(x_l,P[x_1,x_{l-1}])=l-1$,  $d(x_l,P[x_{l+1},x_{p-2}])=p-l-1$ and $x_lx_p, x_{p-2}x_l\in D$.

 Let next $yx_l\notin D$. Similarly as in the case $yx_l\in D$ we deduce that $d(x_l,P[x_{l+1},x_{p-2}])=p-l-1$ and $x_lx_p, x_{p-2}x_l\in D$. (11) is proved.

Now using (9) and (11), we obtain a hamiltonian cycle $H=x_1\ldots x_{l-1}x_{p-1}xx_{l+1}\ldots x_{p-2}x_lx_pyx_1$, which is a contradiction.

\noindent\textbf{Case 5.2}. $l=p-2$. Then $x_px_{p-2}\notin D$ and $d(x_{p-2},\{x_{p-1},x_p\})\leq 2$. By considered case $l\leq p-3$, w.l.o.g. we can assume that the vertex $x$ is adjacent with all vertex of $P[x_1,x_{p-3}]$. Then
$$
N^+(x)=\{x_1,x_2,\ldots ,x_{m-1}, x_{p-1}\} \,\, \hbox{and} \,\, N^-(x)=\{x_{m-1},x_{m},\ldots ,x_{p-3},x_{p-1}, x_{p}\}. \eqno (12)
$$
 This together with $\{x_{p-3},x_{p-1},x_p\}\rightarrow x$ implies that $m\geq 3$ and $xx_2\in D$. Now we divide this case into three subcases.

\noindent\textbf{Subcase 5.2.1}. $yx_{2}\in D$. Assume that $yx_{p-2}\notin D$. Then $d^+(y, P[x_2,x_{p-3}])=m-2$ since $y$ and $x_{p-1}$ are not adjacent. From this and $d^-(x_{p-2},\{x,y,x_p\})=0$ it follows that $x_ix_{p-2}, yx_{i+1}\in D$ for some $i\in [1,p-4]$. Therefore $H=x_1\ldots x_{i}x_{p-2}x_{p-1}x_pyx_{i+1}\ldots x_{p-3}xx_1$, which is a contradiction. So, we can assume that $yx_{p-2}\in D$. Now it is easy to see that $x_1$ and $x_{p-2}$ are not adjacent. Indeed, if $x_1x_{p-2}\in D$, then $H=x_1x_{p-2}x_{p-1}x_{p}yx_2 $ $\ldots x_{p-3}xx_1$; and if $x_{p-2}x_1\in D$, then $H=x_1\ldots x_{p-3}xx_{p-1}x_{p}yx_{p-2}x_1$; which is a contradiction. Since $x_{p-2}$ cannot be inserted into $P[x_2,x_{p-3}]$, by Lemma 2 we have  $d(x_{p-2},P[x_2,x_{p-3}])\leq p-3$. On the other hand,
$$
p+1=d(x_{p-2})=d(x_{p-2},P[x_2,x_{p-3}])+d(x_{p-2},\{x_{p-1},x_{p}\})+a(x_{p-2},y)
$$
implies that $d(x_{p-2},P[x_2,x_{p-3}])=p-3$. Hence, by Lemma 2, $x_{p-2}x_2\in D$ and $x_2\ldots x_{p-3}xx_{p-1}x_pyx_{p-2}$ $x_2$ is a cycle of length $n-1$ which does not contain $x_1$. Therefore $x_1$ is a $T$-vertex. Now we consider the vertex $x_1$. Observe that if $x_1x_i\in D$, $i\in [m, p-2]$, then by (12), $H=x_1x_i \ldots x_p yx_2\ldots x_{i-1}xx_1$; and if $x_1 x_{p-1}\in D$, then $H=x_1x_{p-1}x_pyx_{p-2}x_2 \ldots x_{p-3}xx_1$ a contradiction. Therefore $d^+(x_1,\{x,y,x_m,x_{m+1},$ $\ldots ,x_{p-1}\})=0$ which contradicts that $x$ is a $T$-vertex.

\noindent\textbf{Subcase 5.2.2}. The vertices $x_2$ and $y$ are not adjacent. Then $x_1y, yx_3\in D$ by Claim 1(iii), and by Claim 3 the vertices $x_2$ and $x_p$ also are not adjacent. Observe that if $x_ix\in D$ with $i\in [3,p-1]$, then $x_2x_{i+1}\notin D$ (for otherwise, $H=x_1x_2x_{i+1}\ldots x_{p}yx_{3}\ldots x_{i}xx_1$). From this we have, if $x_2x\notin D$, then $d^-(x,P[x_3,x_{p-1}])=m-1$ and at least $m+2$ vertices are not dominated by $x_2$ since $d^+(x_2, \{y,x,x_1\})=0$, which contradicts that $x_2$ is a $T$-vertex. So, we can assume that                                                                                                                                                                       $x_2x\in D$. Since the vertex $x$ is adjacent with all vertices of $P[x_1,x_{p-3}]$ it follows that $m=3$. Note that $x_2x_4\in D$ by (9), and $x_2, x_3, x_4$ are $T$-vertices. It is easy to see that 
$$
d^+(x_2, \{x_1, x_5,y\})=d^+(x_3, \{x, x_1,x_2\})=d^+(x_4, \{y,x_3, x_1\})=d^-(x_3, \{x,x_4, x_5\})=0.
$$
Therefore $x_3x_5, x_4x_2, x_1x_3\in D$. Since $x_1yx_3x_4x_2xx_1$ (respectively, $x_2x_3yx_5xx_4x_2$) is a cycle of length $n-1=6$, it follows that $x_5$ (respectively, $x_1$) is a $T$-vertex. Now from 
$$d^+(x_5, \{x_2, x_3,x_4\})=d^-(x_1, \{x_2, x_3,x_4\})=0$$ 
we have $x_5x_1\in D$. Therefore $D$ is isomorphic to well-known digraph $D_7$ or is hamiltonian, a contradiction to our assumption. 

\noindent\textbf{Subcase 5.2.3}. $x_2y\in D$ and $yx_2\notin D$. Then by Claim 1(ii) we have $ x_1y\in D$ and there is a vertex $x_k$ with $k\in [3,p-3]$ which is not adjacent with $y$ ( since $m\geq 3$). Then $x_{k-1}y$ and $yx_{k+1}\in D$ by Claim 1(iii). Using Claim 3, we obtain that $x_k$ is not adjacent with $x_1$ and $x_p$. Since  $x_k$ cannot be inserted  into $P[x_{2},x_{k-1}]$  and $P[x_{k+1},x_{p-1}]$, applying Lemma 2 to these paths we obtain that 
 $$
d(x_k,P[x_{2},x_{k-1}])\leq k-1 ,\,\,\, d(x_k,P[x_{k+1},x_{p-1}])\leq p-k,
$$
 $$ p+1\leq d(x_k)=d(x_k,P[x_{2},x_{k-1}])+ d(x_k,P[x_{k+1},x_{p-1}])+a(x_k,x)$$
and $a(x_k,x)=2$ (in other words $xx_k, x_kx\in D$) and each inequality is, in fact, an equality. Hence, by Lemma 2, $x_kx_2, x_{p-1}x_k\in D$. From $xx_k, x_kx\in D$ we obtain that
$$ 
N^+(x)=\{x_1,x_2,\ldots ,x_k,x_{p-1}\} \, \, \hbox{and} \,\, N^-(x)= \{x_k,x_{k+1}, \ldots ,x_{p-3},x_{p-1},x_p\}
$$
and $x_1\ldots x_{k-1}yx_{k+1}\ldots x_{p-1}x_kxx_1$ is a cycle of length $n-1$. Therefore $x_p$ is a $T$-vertex and $k=m-1$. Now we will consider the vertex $x_p$. Then $x_px_i\notin D$ for all $i\in [k,p-1]\cup \{2\}$ (for otherwise, $H=x_1x_2\ldots $ $ x_{i-1}xx_{p-1}x_px_i\ldots x_{p-2}yx_1$ when $i\in [k+1,p-2]$;  and   $H=x_1\ldots x_{i-1}yx_px_{i}\ldots x_{p-1}xx_1$ when $i=2,k,p-1$ which is a contradiction). Thus we have that the vertex $x_p$ does not dominate at least $m+1$ vertices, which is a contradiction since $x_p$ is a $T$-vertex. This contradiction completes the proof of Claim 5. \fbox \\\\ 

By Claim 2 there is a vertex $x_l$, where $l\in [2,p-1]$, which is not adjacent with $x$, and by Claim 1(iii), $x_{l-1}x, xx_{l+1}\in D$.\\

\noindent\textbf{Remark 1}. Let a vertex $x_k$, where $k\in [2,p-1]$ is not adjacent with the vertices $x$ and $y$ (in other words $d(x_k,\{x,y\})=0$). Then $x_px_k, x_kx_1\in D$ and $N^-(x)=N^-(y)$, $N^+(x)=N^+(y)$.

By Claim 1(iii), $x_{k-1}\rightarrow \{x,y\}\rightarrow x_{k+1}$, $x_k$ is a $T$-vertex and $x_k$ cannot be inserted into $P[x_1,x_{k-1}]$ and $P[x_{k+1},x_{p}]$. Using Lemma 2 we obtain that 
$$d(x_k,P[x_1,x_{k-1}])\leq k \, \,\, \hbox{and} \,\,\, d(x_k,P[x_{k+1},x_{p}])\leq p-k+1,$$ 
$$
p+1=d(x_k)=d(x_k,P[x_1,x_{k-1}])+ d(x_k,P[x_{k+1},x_{p}])\leq p+1.
$$
Therefore each inequality is, in fact, an equality. Hence, by Lemma 2, $x_px_k, x_kx_1\in D$.

Now we show that $N^-(x)=N^-(y)$ and $N^+(x)=N^+(y)$. Assume that this is not the case. Let $x_ix\in D$ and $x_iy\notin D$. Then $x_i\notin \{x_{k-1},x_p\}$, and by Claim 1(ii), $yx_{i+1}\in D$. Since $x_kx_1, x_kx_p\in D$, it is not difficult to see that $H=x_1x_2\ldots x_ixx_{k+1}\ldots x_py x_{i+1}\ldots x_kx_1$ when $i<k-1$ and $H=x_1x_2\ldots x_{k-1}yx_{i+1}\ldots x_px_k\ldots x_ixx_1$ when $i>k$ a contradiction. To show that $N^+(x)=N^+(y)$ it suffices to consider the converse digraph of $D$. \fbox \\\\

  \noindent\textbf{Claim 6}. $d^+(x_{p-1},\{x,y\})\leq 1$.

  \noindent\textbf{Proof}. Suppose, on the contrary, that $x_{p-1}x$ and $x_{p-1}y\in D$. Then $l\leq p-2$. Since $D$ is not hamiltonian it follows that if $xx_{i+1}\in D$ or $yx_{i+1}\in D$, then $x_ix_p\notin D$. This together with $d^-(x_p,\{x,y\})=0$ and $d^+(x,P[x_2, x_{p-1}])=m-1$ implies that at least $m+1$ vertices are not dominate $x_p$. Clearly, $x_p$ is not $T$-vertex. We will distinguish three cases according as $x_ly\in D$ or $x_ly\notin D$ and $yx_l\in D$ or $x_l$ and $y$ are not adjacent.

 \noindent\textbf{Case 6.1}. $x_ly\in D$. Then $d^-(x_l, \{x_p,x_{p-1}\})=0$ (for otherwise, if $x_px_l\in D$, then $H=x_1\ldots x_{l-1}xx_{l+1}$ $\ldots x_px_lyx_1$; and if $x_{p-1}x_l\in D$, then $x_1\ldots x_{l-1}xx_{l+1}\ldots x_{p-1}x_lyx_1$ is an $T$-cycle, a contradiction). So, by the above observation we have that $x_p$ and $x_l$ are not adjacent. Since $x_{p-1}x_l\notin D$ and the vertices $x_l$ cannot be inserted into $P[x_1,x_{l-1}]$ and  $P[x_{l+1},x_{p-1}]$, using Lemma 2 we obtain that 

$$d(x_l,P[x_1,x_{l-1}])\leq l \, \,\, \hbox{and} \,\,\, d(x_l,P[x_{l+1},x_{p-1}])\leq p-l-1.$$ 
Therefore
$$p+1=d(x_l)=d(x_l,P[x_1,x_{l-1}])+ d(x_l,P[x_{l+1},x_{p-1}])+a(x_l,y).$$ 
From this we conclude that $yx_l\in D$ and each inequality is, in fact, an equality. Hence, by Lemma 2, $x_lx_1\in D$ and $H=x_1\ldots x_{l-1}xx_{l+1}\ldots x_pyx_lx_1$, which is a contradiction.

\noindent\textbf{Case 6.2}. $x_ly\notin D$ and $yx_l\in D$. Then $x_lx_1\notin D$ (for otherwise, $H=x_1\ldots x_{l-1}xx_{l+1}\ldots x_pyx_lx_1$) and from $d(y)=n-1$ by Claim 1(ii) we have, $yx_{l+1}\in D$. Since $x_l$ cannot be inserted into  $P[x_{l+1},x_{p}]$ and  into $P[x_1,x_{l-1}]$, using Lemma 2 we obtain that
$$
d(x_l,P[x_1,x_{l-1}])= l-1 \, \,\, \hbox{and} \,\,\, d(x_l,P[x_{l+1},x_{p}])=p-l+1,
$$
and $x_px_l\in D$. By Claim 2 there is a vertex $x_k$, where $k\in [2,p-2]$, which is not adjacent with $y$. Then $x_{k-1}y, yx_{k+1}\in D$ (by Claim 1(iii)) and $x_k$ is a $T$-vertex. We can assume that $x_kx\notin D$ (for otherwise, for the vertex $y$ we would have Case 6.1).

First assume that $k\leq l-1$. Then from $x_kx\notin D$ it follows that  $k\leq l-2$. We now will consider the vertex  $x_k$. It is easy to see that $x_kx_p\notin D$ since $D$ is not hamiltonian. Since $x_p$ is not $T$-vertex and $yx_l\in D$ it follows that if $x_px_k\in D$, then $H=x_1\ldots x_{k-1}yx_l\ldots x_px_k\ldots x_{l-1}xx_1$ is a hamiltonian cycle, and if  $x_{p-1}x_k\in D$, then  $x_1\ldots x_{k-1}yx_l\ldots x_{p-1}x_k\ldots x_{l-1}xx_1$ is an $T$-cycle. In each case we have a contradiction. Therefore the vertices $x_k$ and $x_p$ are not adjacent and $x_{p-1}x_k\notin D$. Consequently, since $x_k$ cannot be inserted into $P[x_1,x_{k-1}]$ and $P[x_{k+1},x_{p-1}]$ by Lemma 2 we obtain 
$$d(x_k,P[x_1,x_{k-1}])\leq k \,\,\, \hbox{and}\,\,\, d(x_k,P[x_{k+1},x_{p-1}])\leq p-k-1.$$ 
Therefore
$$
p+1=d(x_k)=d(x_k,P[x_1,x_{k-1}])+ d(x_k,P[x_{k+1},x_{p-1}])+a(x_k,x)\leq k+p-k-1+1=p,
$$ 
which leads to a contradiction since $x_kx\notin D$ ($a(x_k,x)\leq 1$).\\

Second assume that $k\geq l+1$. From $x_ly\notin D$ it follows that $k\geq l+2$. We may assume that $y$ is adjacent with all vertices of $P[x_1,x_{l+1}]$. Then 
$$
\{x_1,x_2, \ldots ,x_{l+1}\}\subseteq N^+(y) \,\,\, \hbox{and}\,\,\, d^-(y,P[x_{l+1},x_{p-1}])= m-1.
$$
Now consider the vertex $x_l$. It is not difficult to see that if $x_iy\in D$, $i\in [l+1,p-1]$, then $x_lx_{i+1}\notin D$ ( for otherwise, $H=x_1\ldots x_lx_{i+1}\ldots x_pxx_l\ldots x_iyx_1$). Therefore since $x_l$ is a $T$-vertex and $d^+(x_l,\{x,y\})=0$, we obtain that $x_l$ does not dominate at least $m+1$ vertices, which is a contradiction and completes the proof of Case 6.2.\\

Let $\{x_{l_1},x_{l_2}, \ldots x_{l_r}\}$ be a set of  vertices which at the same time are not adjacent with $x$ and $y$, where $2\leq l_1< l_2< \cdots <l_r\leq p-1$. Note that (by Claim 1(iii)) for all $i\in [1,r]$ we have $x_{l_i-1}x, xx_{l_i+1}, x_{l_i-1}y$ and $yx_{l_i+1}\in D$.\\

\noindent\textbf{Remark 2}. The set $\{x,y,x_{l_1},x_{l_2}, \ldots ,x_{l_r}\}$ is an independent set of vertices.

Indeed, if $x_{l_i}x_{l_j}\in D$ and $l_i<l_j$, then $H = x_1\ldots x_{l_i}x_{l_j}\ldots x_pxx_{l_i+1}\ldots x_{l_j-1}yx_1$; and if $x_{l_i}x_{l_j}\in D$ and $l_i>l_j$, then by Remark 1, $x_px_{l_i}\in D$ and $H = x_1\ldots x_{l_j-1}yx_{l_i+1}\ldots x_px_{l_i}x_{l_j}\ldots x_{l_i-1}xx_1$. In each case we arrive at a contradiction. \fbox \\\\

\noindent\textbf{Case 6.3}. The vertices $x_l$ and $y$ are not adjacent. We can assume that  for all $j\in [2,p-2]$ the vertices $x_j$ and $x$ are not adjacent if and only if $x_j$ and $y$ are not adjacent. Then by Remarks 1 and 2 for all $i\in [1,r]$ we have 
$$
N^+(x)=N^+(y)=N^+(x_{l_i}) \,\,\ \hbox{and}\,\,\, N^-(x)=N^-(y)=N^-(x_{l_i}),
$$
and  $\{x,y,x_{l_1},x_{l_2}, \ldots ,x_{l_r}\}$ is an independent set of vertices. Not that if $xx_{i+1}\in D$, then $x_ix_p\notin D$ (for otherwise, $H=x_1\ldots x_ix_pxx_{i+1}\ldots x_{p-1}yx_1$). From this and $d^-(x_p,\{x,y\})=0$ it follows that at least $m+1$ vertices are not dominate $x_p$. Therefore, $x_p$ is not $T$-vertex. Similarly, we can show that if $\{x_i,x_{i+1}\}\rightarrow x$ (respectively, $x\rightarrow \{x_j,x_{j+1}\}$), then $x_{i+1}$ (respectively, $x_j$) is not $T$-vertex; and if $xx_i\in D$ and $x_jx\in D$, then $x_{i-1}x_{j+1}\notin D$. The proof of Claim 6 is completed. \fbox\\\\

\noindent\textbf{Claim 7}. $x_{p-1},x\notin D$.

\noindent\textbf{Proof}. Suppose, on the contrary, that $x_{p-1}x\in D$. Then by Claims 5 and 6 we have $x_{p-1}y\notin D$ and $yx_{p-1}\in D$. Hence by Claim 1(ii), $yx_p\in D$. From this and Claim 2 it follows that $m\geq 3$. There are three possibilities: $xx_2\in D$ or $x$ and $x_2$ are not adjacent or $x_2x\in D$.\\  

\noindent\textbf{Case 7.1}. $xx_2\in D$. If $yx_2\in D$ or $y$ and $x_2$ are not adjacent, then for the converse digraph of $D$ we have that Claim 5 or Claim 6 is not true. Thus we can assume that $x_2y\in D$ and $yx_2\notin D$. Then $x_1y\in D$, by Claim 1(ii). Recall that there is a vertex $x_k$ with $k\in [3,p-2]$ (by Claim 2) which is not adjacent with the vertex $y$ and hence by Claim 1(iii), $x_{k-1}y, yx_{k+1}\in D$ and $x_k$ is a $T$-vertex. 

Now we will prove that the vertex $x_k$ is not adjacent with the vertices $x_1$ and $x_p$ and 
$$
x_{p-1}x_k,\, x_{k}x_2,\, x_{k}x,\, xx_k\in D. \eqno (13) $$
 Suppose that this is not the case. If $x_kx_1\in D$, then $H=x_1yx_{k+1}\ldots x_pxx_2\ldots x_kx_1$; if $x_1x_k\in D$, then $H=x_1x_{k}\ldots x_pxx_2\ldots x_{k-1}yx_1$; if $x_kx_p\in D$, then $H=x_1\ldots x_{k}x_pyx_{k+1}\ldots x_{p-1}xx_1$; and finally if $x_px_k\in D$, then $H=x_1\ldots x_{k-1}yx_px_k\ldots x_{p-1}xx_1$. In each case we have a contradiction. Therefore $x_k$ is not adjacent with the vertices $x_1$ and $x_p$. From this it follows that (since $x_k$ is a $T$-vertex)
$$
p+1=d(x_k)=d(x_k,P[x_2,x_{k-1}])+ d(x_k,P[x_{k+1},x_{p-1}])+ a(x_k,x). \eqno (14)
$$
Since the vertex $x_k$ cannot be inserted into $P[x_2,x_{k-1}]$ and $P[x_{k+1},x_{p-1}]$ by Lemma 2 we have, 
  $$d(x_k,P[x_2,x_{k-1}])\leq k-1 \,\,\, \hbox{and} \,\,\, d(x_k,P[x_{k+1},x_{p-1}])\leq p-k.$$

This together with (14) implies that the above inequalities, in fact, are equalities and $a(x,x_k)=2$ (in other words $x_kx, xx_k\in D$). Again using Lemma 2, we obtain that $x_{p-1}x_k,\, x_{k}x_2\in D$. (13) is proved.

From (13) and Claim 2 it follows that $m\geq 4$. By (13), the cycle $x_1\ldots x_{k-1}yx_{k+1}\ldots x_{p-1}x_kxx_1$ (respectively, $x_2\ldots x_{k-1}yx_{k+1}\ldots x_{p}xx_kx_2$) has length $n-1$ and does not contain $x_p$ (respectively, $x_1$). Therefore, $x_p$ and $x_1$ are  $T$-vertices. It is easy to see that 
$$ 
 \hbox{if}\quad yx_i\in D \quad \hbox{with} \quad i\in [2,p-1], \quad \hbox{then}\quad x_{i-1}x_p\notin D \eqno (15)
$$
(otherwise, if $yx_i$ and  $x_{i-1}x_p\in D$, then  $x_1\ldots x_{i-1}x_pyx_{i}\ldots x_{p-1}xx_1$ is a hamiltonian cycle). Note that $x_{k-1}x_p\notin D$ (otherwise if $x_{k-1}x_p\in D$, then by (13), $x_1\ldots x_{k-1}x_pyx_{k+1}\ldots x_{p-1}x_kxx_1$ is a hamiltonian cycle, a contradiction). From (15), $d^+(y,P[x_2,x_{p-1}])=m-2$, $x_{k-1}x_p\notin D$ and $xx_p\notin D$ it follows that at least $m$ vertices are not dominate $x_p$. Consequently, the vertex $y$ is adjacent with all vertices of $P-\{x_k\}$. Hence
$$ 
\{x_1,x_2,\ldots ,x_{k-1}\}\rightarrow y \rightarrow \{x_{k+1},x_{k+2},\ldots ,x_{p}\}, \eqno (16)
$$
and $k-1=p-k=m-1$. From $x_{k-1}x_p\notin D$ and (15), (16) we have
$$
d^-(x_p,P[x_{k-1},x_{p-2}])=0 \,\,\, \hbox{and} \,\,\, \{x_1,x_2,\ldots ,x_{k-2}\} \rightarrow x_p. \eqno (17)
$$ 
 From this and (13) we have that $x_1\ldots x_{k-2}x_pyx_{k+1}\ldots  x_{p-1}x_{k}xx_1$ is a cycle of length $n-1$ which does not contain $x_{k-1}$. This means that $x_{k-1}$ is a $T$-vertex and  $x_{k-1}$ cannot be inserted into $P[x_1,x_{k-2}]$ and $P[x_{k+1},x_{p-1}]x_k$. 

Now we will consider the vertex $x_{k-1}$ and claim that $x_{k-1}$ is not adjacent with the vertices $x_1$ and $x_p$. 
Indeed, if $x_1x_{k-1}\in D$, then by (13), $H=x_1 x_{k-1}\ldots x_pxx_2\ldots x_{k-2}yx_1$; if $x_{k-1}x_1\in D$, then by (17) and (13),  $H=x_1x_py x_{k+1}\ldots x_{p-1}x_{k}xx_2\ldots x_{k-1}x_{1}$; if $x_px_{k-1}\in D$, then by (16), $H=x_1\ldots x_{k-2}yx_px_{k-1} $ $\ldots  x_{p-1}xx_1$; if $x_{k-1}x_p\in D$, then by (13) and (16), $H=x_1\ldots x_{k-1}x_pyx_{k+1}\ldots $ $ x_{p-1}x_{k}xx_1$. In each case we have obtained a contradiction. Therefore $x_{k-1}$ is not adjacent with the vertices $x_1$ and $x_p$. 

Now by Lemma 2 we have 
$$
p+1=d(x_{k-1})=d(x_{k-1}, P[x_2,x_{k-2}])+d(x_{k-1},P[x_{k+1},x_{p-1}]\cup \{x_k\})+a(x_{k-1},\{x,y\})\leq 
$$
$$
p-1+a(x_{k-1},\{x,y\}).
$$
It is possible only if $a(x_{k-1},\{x,y\})=2$ (i.e., $x_{k-1}y$ and $xx_{k-1}\in D$ since $yx_{k-1}\notin D$ and $x_{k-1}x\notin D$). It is not difficult to see that $d^-(x_1,P[x_{k-1},x_{p-1}])=0 $ (otherwise if $x_ix_1\in D$, $i\in [k,p-1]$, then $H= x_1yx_{i+1}\ldots x_pxx_2\ldots x_ix_1$). Hence $x_{k-2}x_1\in D$ and by (13), $H= x_1yx_{k+1}\ldots x_px x_{k-1}x_kx_2 \ldots x_{k-2}x_1$, which is a contradiction. The contradiction completes the proof of Case 7.1 .

\noindent\textbf{Case 7.2}. The vertices $x$ and $x_2$ are not adjacent. Then by Claim 1(iii), $x_1x$ and $xx_3\in D$. By Claim 4 we have that the vertices $x_2$ and $y$ are adjacent. If we consider the converse digraph of $D$, then using Claim 5 we see that $x_2y\in D$ and $yx_2\notin D$. Therefore, by Claim 1(ii), $x_1y\in D$ since $y$ is a $T$-vertex. Now we will consider the vertex $x_2$. Note that $x_2$ also is a $T$-vertex. If $x_px_2\in D$, then $H=x_1yx_px_2\ldots x_{p-1}xx_1$, a contradiction. So, we can assume that $x_px_2\notin D$. By Lemma 2, $d(x_2,P[x_3,x_{p}])\leq p-2$ since $x_2$ cannot be inserted into $P[x_3,x_p]$. From this, since $x$ and $x_2$ are not adjacent, $yx_2\notin D$ and $x_2$ is a $T$-vertex, we obtain that $x_2x_1\in D$. Now it is easy to see that if $yx_i\in D$ with $i\in [4,p]$, then $x_{i-1}x_2\notin D$ (for otherwise, $H=x_1yx_i\ldots x_{p}xx_3\ldots x_{i-1}x_2x_1$). Consequently, from $d^+(y,P[x_4,x_p])=m-1$ and $d^-(x_2,\{x,y\})=0$ it follows that at least $m+1$ vertices are not dominate $x_2$, which is a contradiction. The obtained contradiction completes the proof of Case 7.2 .

\noindent\textbf{Case 7.3}. $x_2x\in D$. Then $x_1x\in D$ by Claim 1(ii). Then from $d^-(x,\{x_1,x_2,x_{p-1},x_p\})=4$ we have $m\geq 4$. It follows that there is a $l\in [3,p-2]$ such that $x_{l-2}x, x_{l-1}x, xx_{l+1}\in D$ and $x_l$ and $x$ are not adjacent by Claim 2. Note that respect to vertices $x_2$ and $y$ the following subcases are possible:  $yx_2\in D$ or $x_2y\in D$  or the vertices $y$ and $x_2$ are not adjacent.

\noindent\textbf{Subcase 7.3.1}. $yx_2\in D$. It is not difficult to see that the vertices $x_1$ and $x_l$ are not adjacent. Indeed, if $x_1x_l\in D$, then $H=x_1x_l\ldots x_{p}yx_2\ldots x_{l-1}xx_1$; and if $x_lx_1\in D$, then $H=x_1xx_{l+1}\ldots x_{p}yx_2\ldots x_{l}x_1$, which is a contradiction. 

We first prove that  
$$
yx_l,\, x_lx_2,\, x_lx_{l-1}, \, x_{l}x_{l-2}\in D \,\, \hbox{and} \, \, x_{l-2}x_{l}\notin D. \eqno (19)
$$
\noindent\textbf{Proof of (19)}. Assume that $x_px_l\in D$. Then $x_ly\notin D$ (for otherwise, if $x_ly\in D$, then $H=x_1\ldots x_{l-1}xx_{l+1}\ldots x_{p}x_lyx_1$). Since $x_1$ and $x_l$ are  not adjacent and $x_l$ cannot be inserted into $P[x_2,x_{l-1}]$ and $P[x_{l+1},x_p]$, using Lemma 2 we see that 

$$
p+1=d(x_{l})=d(x_{l}, P[x_2,x_{l-1}])+d(x_{l},P[x_{l+1},x_{p}])+a(x_{l},y)\leq p+a(x_{l},y).
$$
It follows that $d(x_{l}, P[x_2,x_{l-1}])=l-1$ and $ a(x_{l},y)=1$. Therefore $yx_l\in D$ and $x_lx_2\in D$ by Lemma 2.

Now assume that  $x_px_l\notin D$. Then similarly as before we obtain that $d(x_{l}, P[x_2,x_{l-1}])=l-1$, $d(x_{l},P[x_{l+1},x_{p}])=p-l$ and $a(x_{l},y)=2$ (i.e., $yx_l,x_ly\in D$). By Lemma 2 we have, $x_lx_2\in D$. Now we will consider the path $x_{l+1}x_{l+2}\ldots x_pyx_1\ldots $ $ x_{l-2}x_{l-1}$ and the vertex $x_l$ instead  of $y$. Then using Claims 6 and  5 we obtain that $x_{l}x_{l-1}, x_{l}x_{l-2}\in D$ and $x_{l-2}x_l\notin D$. So indeed(19) satisfied, as desired. \fbox \\\\

W.l.o.g. we can assume that $xx_{l+2}\notin D$ and  $x$ and $x_{l+2}$ are adjacent (because otherwise for the path $x_{l+1}x_{l+2}\ldots x_pyx_1\ldots x_{l-1}$ we would have Case 7.1 or 7.2 which we have already dealt with). Then by Claim 1(ii) we have, $x_{l+1}x, x_{l+2}x\in D$.

Now we consider the vertex $x_1$. If $x_ix\in D$ with $i\in [2,p-1]$, then $x_1x_{i+1}\notin D$ (for otherwise,  $H=x_1x_{i+1}\ldots x_{p}yx_2\ldots x_{i}xx_1$). If $x_1x_{l+1}\in D$, then $H=x_1x_{l+1}\ldots x_{p}yx_lx_2\ldots x_{l-1}xx_1$ by (19). Observe that $x_2\ldots x_{l-1}xx_{l+1}\ldots x_{p}yx_lx_2$ is a cycle of length $n-1$ which does not contain $x_1$. This means that $x_1$ is a $T$-vertex. Now from $d^-(x,P[x_2,x_{p-1}])=m-2$ and $d^+(x_1,\{y,x_{l+1}\})=0$ it follows that the vertex $x$ is adjacent with all vertices of $P-\{x_l\}$ which is not possible since $m\geq 4$, $x_{l+1}x\in D$ and $D$ is not hamiltonian.

\noindent\textbf{Subcase 7.3.2}. $x_2y\in D$. Then by Claims 2 and 1(iii) there is a vertex $x_k$ with $k\in [3,p-2]$ such that $x_{k-1}y, yx_{k+1}\in D$ and $y$ is not adjacent with $x_k$. It is easy to see that $x_p$ and $x_k$ are not adjacent (i.e., $a(x_k,x_p)=0$). Indeed, if $x_kx_p\in D$, then $H=x_1\ldots x_kx_pyx_{k+1}\ldots x_{p-1}xx_1$; and if $x_px_k\in D$, then $H=x_1\ldots x_{k-1}yx_{p}x_k\ldots x_{p-1}x$ $x_1$, which is a contradiction. Now we prove that
$$
x_{p-1}x_k \,\,\, \hbox{and}\,\,\, x_kx\in D. \eqno (20)
$$
\noindent\textbf{Proof of (20)}. Let $x_kx_1\in D$. Then $xx_k\notin D$ (since otherwise if $xx_k\in D$, then $H=x_1\ldots x_{k-1}yx_{k+1}\ldots$ $x_pxx_{k}x_1$) and hence, since $a(x_k,x_p)=0$ and the paths $P[x_1,x_{k-1}]$ and $P[x_{k+1},x_{p-1}]$  cannot be extended with $x_k$ by Lemma 2 we have $d(x_k,P[x_1,x_{k-1}])\leq k$,\, $d(x_{k},P[x_{k+1},x_{p-1}])\leq p-k$ and
$$
p+1=d(x_{k})=d(x_{k}, P[x_1,x_{k-1}])+d(x_{k},P[x_{k+1},x_{p-1}])+a(x_{k},x)=p+1.
$$
Therefore $d(x_k,P[x_1,x_{k-1}])=k$, $d(x_{k},P[x_{k+1},x_{p-1}])=p-k$ and $a(x_k,x)=1$ (i.e., $x_kx\in D $). Now using Lemma 2 we obtain that $x_{p-1}x_k\in D$.

Let now $x_kx_1\notin D$. Then $d(x_{k}, P[x_1,x_{k-1}])\leq k-1$, $a(x_k,x)=2$ (i.e., $x_kx, xx_k\in D$) and $d(x_{k},P[x_{k+1},x_{p-1}])=p-k$. Again using Lemma 2 we obtain that $x_{p-1}x_k\in D$. So indeed (20) is satisfied, as desired. \fbox \\\\

Now we will consider the vertex $x_p$ which is a $T$-vertex since $x_1\ldots x_{k-1}yx_{k+1}\ldots x_{p-1}x_kxx_1$ is a cycle of length $n-1$. If $x_iy\in D$ with $i\in [1,p-2]$, then $x_px_{i+1}\notin D$ (for otherwise,  $H=x_1\ldots x_iyx_{p}x_{i+1}\ldots x_{p-1}xx_1$). Note that $d^-(y,P[x_1,x_{p-2}])=m-1$ and $x_px_{k+1}\notin D$ (if $x_px_{k+1}\in D$, then by (20), $H=x_1\ldots x_{k-1}yx_{p}x_{k+1}\ldots x_{p-1}x_kxx_1$. It follows from the observation above that the vertex $y$ is adjacent with all vertices of $P-\{x_k\}$. Therefore
$$
N^-(y)=\{x_1,x_2, \ldots ,x_{k-1},x_p\}\,\,\, \hbox{and}\,\,\, N^+(y)=\{x_1,x_{k+1},x_{k+2},\ldots ,x_p\}.
$$
Then for the path $x_{k+1}x_{k+2}\ldots x_pxx_1x_2\ldots x_{k-1}$ and for the vertex $y$ by Claims 5 and 6 we have the considered Case 7.1.\\

\noindent\textbf{Subcase 7.3.3}. The vertices $y$ and $x_2$ are not adjacent. Then $x_1y,yx_3\in D$ (by Claim 1(iii)), $x_2$ and $x_p$ are not adjacent (by Claim 3) and $x_2$ is a $T$-vertex. 

Assume that $x_2x_1\in D$. Then $x_ix_2\notin D$ if $xx_{i+1}\in D$, $i\in [3,p-1]$ (for otherwise, $H=x_1x x_{i+1}$ $\ldots x_{p}yx_{3}\ldots $ $ x_{i}x_2x_1$). Now from $d^+(x,P[x_4,x_{p-1}])=m-1$ and $d^-(x_2,\{x,y\})=0$ it follows that $d^-(x_2)\leq m-1$, which is a contradiction. So, we can assume that $x_2x_1\notin D$. Therefore
 
$$
p+1=d(x_{2})=d(x_{2}, P[x_3,x_{p-1}])+d(x_{2},\{x_1,x\})\leq d(x_{2}, P[x_3,x_{p-1}])+2.
$$
Hence $d(x_{2}, P[x_3,x_{p-1}])=p-1$. By Lemma 2,  $ x_2$ can be inserted into path $P[x_3,x_{p-1}]$, a contradiction which completes the proof of Claim 7. \fbox \\

Let us now complete the poof of the theorem. Since $D$ is not hamiltonian from Claim 7 and Remark 2 it follows that for any cycle  $C:=x_1x_2\ldots x_{2m}x_1$ of length $n-1=2m$ if $x\notin V(C)$ then $N^+(x)=N^-(x)= \{x_1,x_3,\ldots ,x_{2m-1}\}$ and $\{x_2,x_4,\ldots ,x_{2m},x\}$ is an independent  set of vertices. Therefore $K^*_{m,m+1}\subseteq D \subseteq [K_m+\overline K_{m+1}]^*$. The proof of the Theorem is complete. \fbox \\\\

\noindent\textbf{Remark 3}. Let $D$ be a digraph with vertex set $V(D)=\{x_1,x_2,x_3,x_4,x_5,x,y\}$ such that $N^+(x_1)=\{x_2,x_4\}$, $N^+(x_2)=\{x,y,x_3,x_5\}$, $N^+(x_3)=N^+(x)=N^+(y)=\{x_1,x_2,x_4,\}$, $N^+(x_4)=\{x,y,x_5\}$ and $N^+(x_5)=\{x,y,x_3\}$. It is easy to check that the vertices $x, y, x_2, x_3$ and $x_4$ are $T$-vertices and the vertices $x_1$ and $x_5$ are not $T$-vertices. Moreover, the digraph $D$ is 2-strong and contains no cycle through  $x, y, x_2, x_3$ and $x_4$.

\end{document}